\documentclass[thmsa,12pt,sw20lart]{article}%
\usepackage{amsmath}
\usepackage{amsfonts}
\usepackage{amssymb}
\usepackage{graphicx}%
\newtheorem{theorem}{Theorem}

\newtheorem{lemma}[theorem]{Lemma}

\newtheorem{remark}[theorem]{Remark}

\newenvironment{proof}[1][Proof]{\noindent\textbf{#1.} }{\ \rule{0.5em}{0.5em}}
\begin{document}

\author{Yu-Chu Lin\ \ and\ \ Dong-Ho Tsai}
\title{On a General Linear Nonlocal Curvature Flow of Convex Plane
Curves\thanks{Mathematics\ Subject Classification:\ 35K05, 35K15.}}
\maketitle

\begin{abstract}
Motivated by Pan-Yang \cite{PY} and Ma-Cheng\ \cite{MC},\ we study a general
linear nonlocal curvature flow for convex closed plane curves and discuss the
short time existence and asymptotic convergence behavior\ of the flow.\ 

Due to the linear structure of the flow, this partial differential equation
problem can be resolved using an ordinary differential equation method,
together\ with the help of representation formula for solutions to a linear
heat equation.\ 

\end{abstract}

\section{Introduction.}

Recently there has been some interest in the \emph{nonlocal} flow of
\emph{convex} closed plane curves\footnote{There are already a lot of nonlocal
curvature flow of hypersurfaces,\ especially those related to the so-called
mean curvature flow.\ We will not mention them here.}. See the papers by
Gage\ \cite{GA1},\ Jiang-Pan\ \cite{JP},\ Pan-Yang \cite{PY}%
,\ Ma-Cheng\ \cite{MC}, Ma-Zhu\ \cite{MZ}\ and Lin-Tsai\ \cite{LT1}.\ All of
the above papers deal with the evolution of a given convex\footnote{Here
"convex" always means "strictly convex".\ A convex closed
plane\ curve\ has\ positive curvature everywhere.\ } simple\ closed
plane\ curve $\gamma_{0}$. The general form of the equation is given by%
\begin{equation}
\left\{
\begin{array}
[c]{l}%
\dfrac{\partial X}{\partial t}\left(  \varphi,t\right)  =\left[  F\left(
k\left(  \varphi,t\right)  \right)  -\lambda\left(  t\right)  \right]
\mathbf{N}_{in}\left(  \varphi,t\right)
\vspace{3mm}%
\\
X\left(  \varphi,0\right)  =X_{0}\left(  \varphi\right)  ,\ \ \ \varphi\in
S^{1},
\end{array}
\right.  \label{main}%
\end{equation}
which is a parabolic initial value problem.\ Here $X_{0}\left(  \varphi
\right)  :S^{1}\rightarrow\gamma_{0}\ $is a smooth parametrization of
$\gamma_{0}$; $k\left(  \varphi,t\right)  \ $is the curvature of the evolving
curve $\gamma_{t}=\gamma\left(  \cdot,t\right)  \ $(parametrized by $X\left(
\varphi,t\right)  $)\ at the point $\varphi;\ $and $\mathbf{N}_{in}\left(
\varphi,t\right)  $ is the inward normal of the curve $\gamma_{t}\ $at time
$t$. As for the speed,$\ F\left(  k\right)  \ $is a given function\ of
curvature\ satisfying the\emph{ parabolic condition} $F^{\prime}\left(
z\right)  >0\ $for all $z\ $in its\ domain\ (usually $\left(  0,\infty\right)
$)$\ $and $\lambda\left(  t\right)  \ $is a function of time, which may depend
on certain global quantities of $\gamma_{t},\ $say its length $L\left(
t\right)  \ $or enclosed area $A\left(  t\right)  ,\ $or others\ (see
(\ref{F2-1})\ and\ (\ref{F4}) below).$\ $In such a case the flow has nonlocal
character.\ Note that if $\lambda\left(  t\right)  $ depends on $\gamma
_{t},\ $then it is not known beforehand. We only know $\lambda\left(
0\right)  .$\ 

For $k$-type nonlocal\ flow,\ the following$\ $
\begin{equation}
F\left(  k\right)  -\lambda\left(  t\right)  =k-\frac{2\pi}{L\left(  t\right)
}\ \ \ \text{(area-preserving,\ gradient flow of\ the IPD\ }L^{2}-4\pi
A\text{)} \label{F1}%
\end{equation}
(IPD means \emph{isoperimetric deficit})\ and
\begin{equation}
F\left(  k\right)  -\lambda\left(  t\right)  =k-\frac{L\left(  t\right)
}{2A\left(  t\right)  }\ \ \ \text{(gradient flow of the IPR\ }\frac{L^{2}%
}{4\pi A}\text{)} \label{F2}%
\end{equation}
(IPR means \emph{isoperimetric ratio})\ and \
\begin{equation}
F\left(  k\right)  -\lambda\left(  t\right)  =k-\frac{1}{2\pi}\int
_{0}^{L\left(  t\right)  }k^{2}ds\ \ \ \text{(length-preserving)} \label{F2-1}%
\end{equation}
have been studied by\ Gage \cite{GA1}, Jiang-Pan\ \cite{JP}%
,\ and\ Ma-Zhu\ \cite{MZ} respectively.\ 

On the other hand, for $1/k$-type\ nonlocal\ flow\ (here the initial curve
$\gamma_{0}\ $must be \emph{convex}, otherwise $k=0$ somewhere), the
following
\begin{equation}
F\left(  k\right)  -\lambda\left(  t\right)  =\frac{1}{L\left(  t\right)
}\int_{0}^{L\left(  t\right)  }\frac{1}{k}ds-\frac{1}{k}%
\ \ \ \text{(area-preserving)} \label{F4}%
\end{equation}
and$\ $
\begin{equation}
F\left(  k\right)  -\lambda\left(  t\right)  =\frac{L\left(  t\right)  }{2\pi
}-\frac{1}{k}\ \ \ \text{(length-preserving,\ dual flow of Gage)} \label{F3}%
\end{equation}
and%
\begin{equation}
F\left(  k\right)  -\lambda\left(  t\right)  =\dfrac{2A\left(  t\right)
}{L\left(  t\right)  }-\dfrac{1}{k}\ \ \ \text{(dual flow of Jiang-Pan)}
\label{F8}%
\end{equation}
have been studied by\ Ma-Cheng\ \cite{MC},\ Pan-Yang \cite{PY},\ and Lin-Tsai
\cite{LT1}\ respectively.\ The good thing about$\ $these$\ 1/k$-type\ flows is
that they produce a \textbf{linear} equation for the radius of curvature
$1/k$, hence it is easier to deal with it. See the discussions below.\ 

When the initial cure $\gamma_{0}$ is convex, the results claimed in each of
the above mentioned papers are roughly more or less the same:\ the flow
preserve the convexity of a given initial curve $\gamma_{0}$ and evolve it (if
no singularity forming in finite time)\ to a round circle in $C^{\infty}$
sense as $t\rightarrow\infty.\ $However,\ unfortunately,\ the $C^{\infty}$
convergence proof in most of the above\ papers are usually omitted due to its
cumbersome details.

We know that (see Gage \cite{GA1}) for a family of time-dependent
smooth\ simple\ closed\ curves $X\left(  \varphi,t\right)  :S^{1}\times
\lbrack0,T)\rightarrow\mathbb{R}^{2}$ with \textsf{time variation}%
\begin{equation}
\frac{\partial X}{\partial t}\left(  \varphi,t\right)  =W\left(
\varphi,t\right)  \in\mathbb{R}^{2}, \label{W}%
\end{equation}
its length $L\left(  t\right)  $ and enclosed area $A\left(  t\right)
\ $evolve according to the following:
\begin{equation}
\frac{dL}{dt}\left(  t\right)  =-\int_{\gamma_{t}}\left\langle W,\ k\mathbf{N}%
_{in}\right\rangle ds,\ \ \text{\ \ \ }\frac{dA}{dt}\left(  t\right)
=-\int_{\gamma_{t}}\left\langle W,\ \mathbf{N}_{in}\right\rangle ds,
\label{LA-evo}%
\end{equation}
where $s$ is the arc\ length parameter of the curve $X\left(  \varphi
,t\right)  .\ $As for the curvature $k\left(  \varphi,t\right)  $,\ following
computations similar to those in Gage-Hamilton\ \cite{GH}, we can obtain the
equation%
\begin{equation}
\frac{\partial k}{\partial t}\left(  \varphi,t\right)  =\left\langle
\frac{\partial^{2}W}{\partial s^{2}},\ \mathbf{N}_{in}\right\rangle
-2k\left\langle \frac{\partial W}{\partial s},\ \mathbf{T}\right\rangle
\label{dkdt}%
\end{equation}
where $\mathbf{T}=\mathbf{T}\left(  \varphi,t\right)  $ is the unit tangent
vector of$\ X\left(  \varphi,t\right)  .\ $

\section{A general linear nonlocal curvature flow;\ short time existence of a
solution.}

Regarding the short time existence of a solution to the nonlocal flow
(\ref{main}), for $k$-type nonlocal\ flow in (\ref{F1}),\ (\ref{F2}) or
(\ref{F2-1}), one can use the same method as in\ Section
2\ of\ Gage-Hamilton\ \cite{GH} to show that for any smooth\ convex closed
curve$\ $(or simple closed curve) $\gamma_{0},\ $the flows in (\ref{F1}%
),\ (\ref{F2}), (\ref{F2-1}) all have a smooth solution defined on
$S^{1}\times\lbrack0,T)\ $for short time $T>0.\ $If $\gamma_{0}$ is convex, by
continuity the convexity is also preserved for short time.{}

Let $H\left(  p,q\right)  :\left(  0,\infty\right)  \times\left(
0,\infty\right)  \rightarrow\mathbb{R}$ be a given\ (but can be arbitrary)
smooth function of two variables.\ The goal of this paper is to first use an
ODE method\ to explain that for a general $1/k$-type nonlocal flow of the form
(with $\gamma_{0}$ convex, parametrized by $X_{0}$)%
\begin{equation}
\left\{
\begin{array}
[c]{l}%
\dfrac{\partial X}{\partial t}\left(  \varphi,t\right)  =\left(  H\left(
L\left(  t\right)  ,A\left(  t\right)  \right)  -\dfrac{1}{k\left(
\varphi,t\right)  }\right)  \mathbf{N}_{in}\left(  \varphi,t\right)
\vspace{3mm}%
\\
X\left(  \varphi,0\right)  =X_{0}\left(  \varphi\right)  ,\ \ \ \varphi\in
S^{1}.
\end{array}
\right.  \label{main-1}%
\end{equation}
it has a smooth convex solution for short time $[0,T).\ $After that, using the
linear structure of $1/k\ $(or the support function$\ u$) together\ with the
help of the representation formula for a linear heat equation, we can solve
$1/k$ (or the support function $u$)\ explicitly and give a precise description
of the evolving curve $\gamma_{t}\ $and then discuss its possible\ convergence behavior.\ 

If (\ref{main-1}) has a convex\ solution defined on $S^{1}\times\lbrack0,T)$
for some $T>0,\ $then using normal (or tangent)\ angle $\theta\in S^{1}\ $of
the curve\ as parameter (as done in Gage-Hamilton\ \cite{GH} and many
others),\ the curvature evolution equation (\ref{dkdt})\ becomes%
\begin{align}
\frac{\partial k}{\partial t}\left(  \theta,t\right)   &  =\left\langle
\frac{\partial^{2}W}{\partial s^{2}},\ \mathbf{N}_{in}\right\rangle
-2k\left\langle \frac{\partial W}{\partial s},\ \mathbf{T}\right\rangle
-\left\langle \frac{\partial W}{\partial s},\ \mathbf{N}_{in}\right\rangle
\frac{\partial k}{\partial\theta}\nonumber\\
&  =k^{2}\left(  \theta,t\right)  \left(  -\frac{1}{k\left(  \theta,t\right)
}\right)  _{\theta\theta}+k^{2}\left(  \theta,t\right)  \left(  H\left(
L\left(  t\right)  ,A\left(  t\right)  \right)  -\frac{1}{k\left(
\theta,t\right)  }\right)  \label{kk}%
\end{align}
where we have used $W=\left(  H\left(  L,A\right)  -1/k\right)  \mathbf{N}%
_{in}\ $and the relation $\partial/\partial s=k\partial/\partial\theta\ $in
(\ref{kk}).$\ $The$\ $attractive feature now is that\ the radius of curvature
$1/k\ $satisfies a \textbf{linear} equation on$\ S^{1}\times\lbrack0,T)$
\begin{equation}
\frac{\partial}{\partial t}\frac{1}{k\left(  \theta,t\right)  }=\left(
\frac{1}{k\left(  \theta,t\right)  }\right)  _{\theta\theta}+\frac{1}{k\left(
\theta,t\right)  }-H\left(  L\left(  t\right)  ,A\left(  t\right)  \right)  .
\label{1/k}%
\end{equation}
By (\ref{1/k}), the evolution of $L\left(  t\right)  \ $is given by%
\begin{equation}
\frac{dL}{dt}\left(  t\right)  =L\left(  t\right)  -2\pi H\left(  L\left(
t\right)  ,A\left(  t\right)  \right)  ,\ \ \ \text{where\ }L\left(  t\right)
=\int_{0}^{2\pi}\frac{1}{k\left(  \theta,t\right)  }d\theta\label{LODE}%
\end{equation}
(which is not self-contained because of the term $A\left(  t\right)  $)\ and
then%
\begin{equation}
\frac{\partial}{\partial t}\left(  \frac{1}{k\left(  \theta,t\right)  }%
-\frac{L\left(  t\right)  }{2\pi}\right)  =\left(  \frac{1}{k\left(
\theta,t\right)  }-\frac{L\left(  t\right)  }{2\pi}\right)  _{\theta\theta
}+\left(  \frac{1}{k\left(  \theta,t\right)  }-\frac{L\left(  t\right)  }%
{2\pi}\right)  , \label{kL}%
\end{equation}
where the term $H\left(  L\left(  t\right)  ,A\left(  t\right)  \right)  $ has
disappeared !\ Thus the quantity $V\left(  \theta,t\right)  =\left[
1/k\left(  \theta,t\right)  -L\left(  t\right)  /2\pi\right]  e^{-t}$
satisfies a linear heat equation $V_{t}\left(  \theta,t\right)  =V_{\theta
\theta}\left(  \theta,t\right)  $on $S^{1}\times\lbrack0,T)\ $with$\ V\left(
\theta,0\right)  =1/k\left(  \theta,0\right)  -L\left(  0\right)  /2\pi.\ $By
the representation formula for solution to a heat equation, the
curvature$\ k\left(  \theta,t\right)  \ $satisfies the simple-looking
\emph{integral equation}:%
\begin{equation}
\frac{1}{k\left(  \theta,t\right)  }-\frac{1}{2\pi}\int_{0}^{2\pi}\frac
{1}{k\left(  \theta,t\right)  }d\theta=e^{t}\int_{-\infty}^{\infty}\frac
{1}{2\sqrt{\pi t}}e^{-\frac{\left(  \theta-\xi\right)  ^{2}}{4t}}\left(
\frac{1}{k\left(  \xi,0\right)  }-\frac{1}{2\pi}\int_{0}^{2\pi}\frac
{1}{k\left(  \sigma,0\right)  }d\sigma\right)  d\xi, \label{int-eq}%
\end{equation}
where the right hand side of the above$\ $is known. Note also that the
function$\ 1/k\left(  \theta,t\right)  -L\left(  t\right)  /2\pi$ is,\ if not
identically equal to $0$, somewhere positive and somewhere negative since its
integral over $S^{1}$ is zero.

Unfortunately, (\ref{int-eq})\ is not enough for us to solve$\ k\left(
\theta,t\right)  $ uniquely\ (if $1/k\left(  \theta,t\right)  \ $satisfies
(\ref{int-eq}), so does $1/k\left(  \theta,t\right)  \ $plus any function of
time).\ This is not surprising because we have not used any information given
by the nonlocal term $H\left(  L,A\right)  \ $so far.\ Moreover, even if we
find one $k\left(  \theta,t\right)  \ $satisfying (\ref{int-eq}),\ we do not
know whether it satisfies the evolution equation (\ref{1/k}) or not.\ To
overcome this, we need to explore other geometric\ quantities.\ 

The \emph{support function} $u\left(  \theta,t\right)  $ of a convex closed
curve $\gamma_{t}$ is, by definition, given by%
\begin{equation}
u\left(  \theta,t\right)  =\left\langle X\left(  \theta,t\right)  ,\ \left(
\cos\theta,\sin\theta\right)  \right\rangle ,\ \ \ \theta\in S^{1}%
\end{equation}
where $X\left(  \theta,t\right)  $ is the position vector of the unique point
on $\gamma_{t}$ with outward normal $\mathbf{N}_{out}\ $equal to $\left(
\cos\theta,\sin\theta\right)  .$ Using $u\left(  \theta,t\right)  ,\ $one can
express $L\left(  t\right)  ,\ A\left(  t\right)  \ $and curvature $k\left(
\theta,t\right)  \ $of $\gamma_{t}\ $as (see the book by Schneider \cite{S}
for details)%
\begin{equation}
k\left(  \theta,t\right)  =\frac{1}{u_{\theta\theta}\left(  \theta,t\right)
+u\left(  \theta,t\right)  }%
\end{equation}
and%
\begin{equation}
L\left(  t\right)  =\int_{0}^{2\pi}u\left(  \theta,t\right)  d\theta
,\ \ \ A\left(  t\right)  =\frac{1}{2}\int_{0}^{2\pi}u\left(  \theta,t\right)
\left(  u_{\theta\theta}\left(  \theta,t\right)  +u\left(  \theta,t\right)
\right)  d\theta. \label{LA}%
\end{equation}
Under the nonlocal flow\ (\ref{main-1}), the evolution of $u\left(
\theta,t\right)  $ on $S^{1}\times\lbrack0,T)$ is
\begin{align}
\frac{\partial u}{\partial t}\left(  \theta,t\right)   &  =u_{\theta\theta
}\left(  \theta,t\right)  +u\left(  \theta,t\right)  -H\left(  L\left(
t\right)  ,A\left(  t\right)  \right) \nonumber\\
&  =u_{\theta\theta}\left(  \theta,t\right)  +u\left(  \theta,t\right)
-H\left(  \int_{0}^{2\pi}u\left(  \theta,t\right)  d\theta,\ \frac{1}{2}%
\int_{0}^{2\pi}u\left(  \theta,t\right)  \left(  u_{\theta\theta}\left(
\theta,t\right)  +u\left(  \theta,t\right)  \right)  d\theta\right)  .
\label{dudt}%
\end{align}
Again, we have$\ $%
\begin{equation}
\frac{\partial}{\partial t}\left(  u\left(  \theta,t\right)  -\frac{L\left(
t\right)  }{2\pi}\right)  =\left(  u\left(  \theta,t\right)  -\frac{L\left(
t\right)  }{2\pi}\right)  _{\theta\theta}+\left(  u\left(  \theta,t\right)
-\frac{L\left(  t\right)  }{2\pi}\right)  \label{u-L}%
\end{equation}
and similar to\ (\ref{int-eq})\ we conclude
\begin{equation}
u\left(  \theta,t\right)  -\frac{L\left(  t\right)  }{2\pi}=e^{t}\int
_{-\infty}^{\infty}\frac{1}{2\sqrt{\pi t}}e^{-\frac{\left(  \theta-\xi\right)
^{2}}{4t}}\left(  u\left(  \xi,0\right)  -\frac{1}{2\pi}\int_{0}^{2\pi
}u\left(  \sigma,0\right)  d\sigma\right)  d\xi, \label{uL}%
\end{equation}
where the the right hand side of is known.\ 

Denote the right hand side of (\ref{uL})\ as $B\left(  \theta,t\right)  $. By
(\ref{uL})\ and (\ref{LA}), we can express the equation (\ref{LODE})\ for
$L\left(  t\right)  $ as%
\begin{align}
\frac{dL}{dt}\left(  t\right)   &  =L\left(  t\right)  -2\pi H\left(  L\left(
t\right)  ,\ \frac{1}{2}\int_{0}^{2\pi}\left(  B\left(  \theta,t\right)
+\frac{L\left(  t\right)  }{2\pi}\right)  \left(  B_{\theta\theta}\left(
\theta,t\right)  +B\left(  \theta,t\right)  +\frac{L\left(  t\right)  }{2\pi
}\right)  d\theta\right) \nonumber\\
&  =L\left(  t\right)  -2\pi H\left(  L\left(  t\right)  ,\ \frac{L^{2}\left(
t\right)  }{4\pi}+D\left(  t\right)  L\left(  t\right)  +E\left(  t\right)
\right)  , \label{L-ode}%
\end{align}
where
\begin{equation}
D\left(  t\right)  :=\frac{1}{2\pi}\int_{0}^{2\pi}B\left(  \theta,t\right)
d\theta,\ \ \ E\left(  t\right)  :=\frac{1}{2}\int_{0}^{2\pi}B\left(
\theta,t\right)  \left(  B_{\theta\theta}\left(  \theta,t\right)  +B\left(
\theta,t\right)  \right)  d\theta\label{DE}%
\end{equation}
are both known functions of time, depending only on the initial data. One can
simplify $D\left(  t\right)  $ and $E\left(  t\right)  $ further.\ Note that%
\[
B\left(  \theta,t\right)  =\frac{e^{t}}{2\sqrt{\pi t}}\int_{-\infty}^{\infty
}e^{-\frac{\xi^{2}}{4t}}\left(  u\left(  \theta-\xi,0\right)  -\frac{L\left(
0\right)  }{2\pi}\right)  d\xi=\frac{e^{t}}{2\sqrt{\pi t}}\int_{-\infty
}^{\infty}e^{-\frac{\xi^{2}}{4t}}u\left(  \theta-\xi,0\right)  d\xi
-\frac{L\left(  0\right)  }{2\pi}e^{t}%
\]
and so
\begin{equation}
D\left(  t\right)  =\frac{1}{2\pi}\int_{0}^{2\pi}B\left(  \theta,t\right)
d\theta=\frac{1}{2\pi}\frac{e^{t}}{2\sqrt{\pi t}}\int_{-\infty}^{\infty
}e^{-\frac{\xi^{2}}{4t}}\left(  \int_{0}^{2\pi}\left(  u\left(  \theta
-\xi,0\right)  -\frac{L\left(  0\right)  }{2\pi}\right)  d\theta\right)
d\xi=0. \label{D0}%
\end{equation}
With (\ref{D0}),\ we find that%
\begin{equation}
E\left(  t\right)  =\frac{1}{2}\int_{0}^{2\pi}B\left(  \theta,t\right)
\left(  B_{\theta\theta}\left(  \theta,t\right)  +B\left(  \theta,t\right)
\right)  d\theta=-\frac{1}{4\pi}\left(  L^{2}\left(  t\right)  -4\pi A\left(
t\right)  \right)  \leq0. \label{E}%
\end{equation}
Also by%
\begin{equation}
B_{\theta\theta}\left(  \theta,t\right)  +B\left(  \theta,t\right)
=\frac{e^{t}}{2\sqrt{\pi t}}\int_{-\infty}^{\infty}e^{-\frac{\xi^{2}}{4t}%
}\left(  u_{\theta\theta}\left(  \theta-\xi,0\right)  +u\left(  \theta
-\xi,0\right)  \right)  d\xi-\frac{L\left(  0\right)  }{2\pi}e^{t} \label{BB}%
\end{equation}
we get$\ $%
\[
E\left(  t\right)  =\frac{1}{2}%
{\displaystyle\int_{0}^{2\pi}}
\left\{
\begin{array}
[c]{l}%
\left(  \frac{e^{t}}{2\sqrt{\pi t}}\int_{-\infty}^{\infty}e^{-\frac{\xi^{2}%
}{4t}}u\left(  \theta-\xi,0\right)  d\xi-\frac{L\left(  0\right)  }{2\pi}%
e^{t}\right)
\vspace{3mm}%
\\
\times\left(  \frac{e^{t}}{2\sqrt{\pi t}}\int_{-\infty}^{\infty}e^{-\frac
{\xi^{2}}{4t}}\left[  u_{\theta\theta}\left(  \theta-\xi,0\right)  +u\left(
\theta-\xi,0\right)  \right]  d\xi\right)
\end{array}
\right\}  d\theta=E_{1}\left(  t\right)  -E_{2}\left(  t\right)  ,
\]
where%
\[
E_{1}\left(  t\right)  =\frac{1}{2}\int_{0}^{2\pi}\left\{
\begin{array}
[c]{l}%
\left(  \frac{e^{t}}{2\sqrt{\pi t}}\int_{-\infty}^{\infty}e^{-\frac{\xi^{2}%
}{4t}}u\left(  \theta-\xi,0\right)  d\xi\right) \\
\left(  \frac{e^{t}}{2\sqrt{\pi t}}\int_{-\infty}^{\infty}e^{-\frac{\xi^{2}%
}{4t}}\left[  u_{\theta\theta}\left(  \theta-\xi,0\right)  +u\left(
\theta-\xi,0\right)  \right]  d\xi\right)
\end{array}
\right\}  d\theta
\]
and%
\[
E_{2}\left(  t\right)  =-\frac{1}{2}\int_{0}^{2\pi}\frac{L\left(  0\right)
}{2\pi}e^{t}\left(  \frac{e^{t}}{2\sqrt{\pi t}}\int_{-\infty}^{\infty
}e^{-\frac{\xi^{2}}{4t}}\left[  u_{\theta\theta}\left(  \theta-\xi,0\right)
+u\left(  \theta-\xi,0\right)  \right]  d\xi\right)  d\theta=-\frac
{L^{2}\left(  0\right)  }{4\pi}e^{2t}.
\]
Thus the ODE (\ref{L-ode})\ becomes%
\begin{equation}
\left\{
\begin{array}
[c]{l}%
\dfrac{dL}{dt}\left(  t\right)  =L\left(  t\right)  -2\pi H\left(  L\left(
t\right)  ,\ \dfrac{L^{2}\left(  t\right)  }{4\pi}+E_{1}\left(  t\right)
-\dfrac{L^{2}\left(  0\right)  }{4\pi}e^{2t}\right)
\vspace{3mm}%
\\
L\left(  0\right)  =\text{\ length of }\gamma_{0},\ \ \ \ \ H\left(
p,q\right)  :\left(  0,\infty\right)  \times\left(  0,\infty\right)
\rightarrow\mathbb{R}.
\end{array}
\right.  \label{L-ODE}%
\end{equation}
The ODE (\ref{L-ODE}) for $L\left(  t\right)  $ is now \emph{self-contained}
and by standard ODE theory\ we can conclude the$\ $following immediately:\ 

\begin{lemma}
For the ODE$\ $in (\ref{L-ODE}) with initial condition $L\left(  0\right)
>0,$ there is a unique\ solution $L\left(  t\right)  \ $defined on interval
$[0,T_{\ast})\ $for some $T_{\ast}>0.\ $
\end{lemma}

Now if we want to solve the initial value problem\ (\ref{main-1}), we can go
in the reverse direction and first use the ODE (\ref{L-ODE}) to obtain
$L\left(  t\right)  $ on $[0,T_{\ast}).\ $Then we use (\ref{uL}) to obtain the
support function $u\left(  \theta,t\right)  .\ $As long as the inequality
$0<u_{\theta\theta}\left(  \theta,t\right)  +u\left(  \theta,t\right)
<\infty\ $is satisfied, we can use it to construct a convex closed
smooth\ curve $\gamma_{t}.\ $Finally we claim that this family of curves
$\gamma_{t}$ is indeed a solution to the nonlocal curvature flow
(\ref{main-1}) by checking that its radius of curvature $1/k\ $or support
function $u$ satisfies the right equation.\ 

Let%
\begin{align}
\tilde{u}\left(  \theta,t\right)   &  :=B\left(  \theta,t\right)
+\frac{L\left(  t\right)  }{2\pi}\nonumber\\
&  =\frac{e^{t}}{2\sqrt{\pi t}}\int_{-\infty}^{\infty}e^{-\frac{\xi^{2}}{4t}%
}u\left(  \theta-\xi,0\right)  d\xi+\frac{L\left(  t\right)  }{2\pi}%
-\frac{L\left(  0\right)  }{2\pi}e^{t},\ \ \ \left(  \theta,t\right)  \in
S^{1}\times\lbrack0,T_{\ast}) \label{u-1}%
\end{align}
where $L\left(  t\right)  $ is from the solution\ in\ ODE (\ref{L-ODE}) with
initial condition\ $L\left(  0\right)  =L\left(  \gamma_{0}\right)  $ and
$u\left(  \cdot,0\right)  \ $is the support function of $\gamma_{0}.\ $The
above $\tilde{u}\left(  \theta,t\right)  $ satisfies$\ \tilde{u}\left(
\theta,0\right)  =u\left(  \theta,0\right)  \ $for all\ $\theta\ $and%
\begin{align}
&  \tilde{u}_{\theta\theta}\left(  \theta,t\right)  +\tilde{u}\left(
\theta,t\right) \nonumber\\
&  =\frac{e^{t}}{2\sqrt{\pi t}}\int_{-\infty}^{\infty}e^{-\frac{\xi^{2}}{4t}%
}\left(  u_{\theta\theta}\left(  \theta-\xi,0\right)  +u\left(  \theta
-\xi,0\right)  \right)  d\xi+\frac{L\left(  t\right)  }{2\pi}-\frac{L\left(
0\right)  }{2\pi}e^{t} \label{good}%
\end{align}
with
\[
\tilde{u}_{\theta\theta}\left(  \theta,0\right)  +\tilde{u}\left(
\theta,0\right)  =u_{\theta\theta}\left(  \theta,0\right)  +u\left(
\theta,0\right)  >0\ \text{\ \ for\ all\ \ \ }\theta\in S^{1}.
\]
As $S^{1}\ $is compact and $\tilde{u}\left(  \theta,t\right)  \ $is continuous
on $S^{1}\times\lbrack0,T_{\ast}),\ $by choosing $T_{\ast}\ $smaller if
necessary,\ we have%
\begin{equation}
\tilde{u}_{\theta\theta}\left(  \theta,t\right)  +\tilde{u}\left(
\theta,t\right)  >0\ \ \ \text{for all\ \ \ }\left(  \theta,t\right)  \in
S^{1}\times\lbrack0,T_{\ast}). \label{u-pos}%
\end{equation}
Note that we always have $\tilde{u}_{\theta\theta}\left(  \theta,t\right)
+\tilde{u}\left(  \theta,t\right)  <\infty\ $as long as time is finite.$\ $%
Therefore for each time $t\in\lbrack0,T_{\ast})\ $there exists a smooth convex
closed curve $\tilde{\gamma}_{t}\ $with support function equal to the above
$\tilde{u}\left(  \theta,t\right)  \ $(see Schneider \cite{S}\ or Lemma 2 in
\cite{LT3}). In fact, the position vector $\tilde{P}\left(  \theta,t\right)
\in\mathbb{R}^{2}$ of $\tilde{\gamma}_{t}$ is given by%
\begin{align}
&  \tilde{P}\left(  \theta,t\right)  =\tilde{u}\left(  \theta,t\right)
\left(  \cos\theta,\sin\theta\right)  +\tilde{u}_{\theta}\left(
\theta,t\right)  \left(  -\sin\theta,\cos\theta\right)  ,\ \ \ \ \ \tilde
{u}\left(  \theta,t\right)  =B\left(  \theta,t\right)  +\frac{L\left(
t\right)  }{2\pi}\nonumber\\
&  =\left\{
\begin{array}
[c]{l}%
e^{t}\int_{-\infty}^{\infty}\frac{1}{2\sqrt{\pi t}}e^{-\frac{\xi^{2}}{4t}%
}\left[  u\left(  \theta-\xi,0\right)  \left(  \cos\theta,\sin\theta\right)
+u_{\theta}\left(  \theta-\xi,0\right)  \left(  -\sin\theta,\cos\theta\right)
\right]  d\xi%
\vspace{3mm}%
\\
+\left(  \frac{L\left(  t\right)  }{2\pi}-\frac{L\left(  0\right)  }{2\pi
}\right)  \left(  \cos\theta,\sin\theta\right)
\end{array}
\right.  \label{position}%
\end{align}
By (\ref{u-1}), the length $\tilde{L}\left(  t\right)  \ $of $\tilde{\gamma
}_{t}$ is given by%
\[
\tilde{L}\left(  t\right)  =\int_{0}^{2\pi}\tilde{u}\left(  \theta,t\right)
d\theta=\int_{0}^{2\pi}\left(  B\left(  \theta,t\right)  +\frac{L\left(
t\right)  }{2\pi}\right)  d\theta=L\left(  t\right)  \ \ \ \text{for
all\ \ \ }t\in\lbrack0,T_{\ast})
\]
due to (\ref{D0}). Also by (\ref{DE}) the enclosed area $\tilde{A}\left(
t\right)  \ $of $\tilde{\gamma}_{t}$ is given by
\begin{align*}
\tilde{A}\left(  t\right)   &  =\frac{1}{2}\int_{0}^{2\pi}\tilde{u}\left(
\theta,t\right)  \left(  \tilde{u}_{\theta\theta}\left(  \theta,t\right)
+\tilde{u}\left(  \theta,t\right)  \right)  d\theta\\
&  =\frac{1}{2}\int_{0}^{2\pi}\left(  B\left(  \theta,t\right)  +\frac
{L\left(  t\right)  }{2\pi}\right)  \left(  B_{\theta\theta}\left(
\theta,t\right)  +B\left(  \theta,t\right)  +\frac{L\left(  t\right)  }{2\pi
}\right)  d\theta\\
&  =\frac{L^{2}\left(  t\right)  }{4\pi}+E_{1}\left(  t\right)  -\frac
{L^{2}\left(  0\right)  }{4\pi}e^{2t}\ \ \ \text{for all\ \ \ }t\in
\lbrack0,T_{\ast}).
\end{align*}
Now the\ support function of $\tilde{\gamma}_{t}$ has the evolution%
\begin{align*}
\frac{\partial\tilde{u}}{\partial t}\left(  \theta,t\right)   &
=\frac{\partial}{\partial t}\left(  B\left(  \theta,t\right)  +\frac{L\left(
t\right)  }{2\pi}\right)  =B_{\theta\theta}\left(  \theta,t\right)  +B\left(
\theta,t\right)  +\frac{\tilde{L}\left(  t\right)  }{2\pi}-H\left(  \tilde
{L}\left(  t\right)  ,\tilde{A}\left(  t\right)  \right) \\
&  =\tilde{u}_{\theta\theta}\left(  \theta,t\right)  +\tilde{u}\left(
\theta,t\right)  -H\left(  \tilde{L}\left(  t\right)  ,\tilde{A}\left(
t\right)  \right)  .
\end{align*}
Finally the curvature$\ \tilde{k}\left(  \theta,t\right)  >0\ $of\ $\tilde
{\gamma}_{t}\ $is given by $\ $%
\begin{equation}
\frac{1}{\tilde{k}\left(  \theta,t\right)  }=\tilde{u}_{\theta\theta}\left(
\theta,t\right)  +\tilde{u}\left(  \theta,t\right)
\end{equation}
and its evolution is
\begin{align}
\frac{\partial}{\partial t}\frac{1}{\tilde{k}\left(  \theta,t\right)  }  &
=\left(  \frac{\partial B}{\partial t}\right)  _{\theta\theta}\left(
\theta,t\right)  +\left(  \frac{\partial B}{\partial t}\right)  \left(
\theta,t\right)  +\frac{\tilde{L}\left(  t\right)  }{2\pi}-H\left(  \tilde
{L}\left(  t\right)  ,\tilde{A}\left(  t\right)  \right) \nonumber\\
&  =\left(  \frac{1}{\tilde{k}\left(  \theta,t\right)  }\right)
_{\theta\theta}+\frac{1}{\tilde{k}\left(  \theta,t\right)  }-H\left(
\tilde{L}\left(  t\right)  ,\tilde{A}\left(  t\right)  \right)  ,\ \ \ \left(
\theta,t\right)  \in S^{1}\times\lbrack0,T_{\ast}) \label{k-evo}%
\end{align}
which\ has the right form and is same\ as the original (\ref{1/k}).

Now we can follow the same argument as in the proof of Theorem 4.1.4 (it says
that the flow equation is equivalent\ to the curvature equation)\ of
Gage-Hamilton\ \cite{GH} to conclude that the nonlocal flow (\ref{main-1}%
)\ has a solution for short time, defined on $S^{1}\times\lbrack0,T_{\ast}%
)\ $for some small $T_{\ast}>0.$

We summarize the following:

\begin{theorem}
(\textbf{short time existence})\ For any smooth function$\ H\left(
p,q\right)  :\left(  0,\infty\right)  \times\left(  0,\infty\right)
\rightarrow\mathbb{R},$ the$\ $nonlocal flow (\ref{main-1}) has a smooth
convex solution defined on $S^{1}\times\lbrack0,T_{\ast})\ $for some time
$T_{\ast}>0.\ $Moreover,\ during this time interval $[0,T_{\ast}),\ $the
parametrization of $\gamma_{t}$ using its outward normal angle $\theta\ $is
given by%
\begin{equation}
P\left(  \theta,t\right)  =\left\{
\begin{array}
[c]{l}%
e^{t}\int_{-\infty}^{\infty}\frac{1}{2\sqrt{\pi t}}e^{-\frac{\xi^{2}}{4t}%
}\left[  u\left(  \theta-\xi,0\right)  \left(  \cos\theta,\sin\theta\right)
+u_{\theta}\left(  \theta-\xi,0\right)  \left(  -\sin\theta,\cos\theta\right)
\right]  d\xi%
\vspace{3mm}%
\\
+\left(  \frac{L\left(  t\right)  }{2\pi}-\frac{L\left(  0\right)  }{2\pi
}\right)  \left(  \cos\theta,\sin\theta\right)  .
\end{array}
\right.
\end{equation}

\end{theorem}

To end this section we point out that under the flow (\ref{main-1}) the
isoperimetric deficit $L^{2}\left(  t\right)  -4\pi A\left(  t\right)  $ is
always decreasing.\ If we compute%
\begin{align*}
&  \frac{d}{dt}\left(  L^{2}\left(  t\right)  -4\pi A\left(  t\right)
\right)  \\
&  =2L\left(  t\right)  \left(  L\left(  t\right)  -2\pi H\left(  L\left(
t\right)  ,A\left(  t\right)  \right)  \right)  -4\pi\left(  \int_{\gamma_{t}%
}\frac{1}{k}ds-H\left(  L\left(  t\right)  ,A\left(  t\right)  \right)
L\left(  t\right)  \right)  \\
&  =2L^{2}\left(  t\right)  -4\pi\int_{\gamma_{t}}\frac{1}{k}ds
\end{align*}
we see that there is \textbf{no} $H\left(  L\left(  t\right)  ,A\left(
t\right)  \right)  \ $term in the time derivative.\ Now by Green-Osher's
inequality\ (also see Pan-Yang \cite{PY})
\begin{equation}
\int_{\gamma_{t}}\frac{1}{k}ds\geq\frac{L^{2}\left(  t\right)  -2\pi A\left(
t\right)  }{\pi}\label{GO1}%
\end{equation}
the above becomes
\[
\frac{d}{dt}\left(  L^{2}\left(  t\right)  -4\pi A\left(  t\right)  \right)
\leq2L^{2}\left(  t\right)  -4\left(  L^{2}\left(  t\right)  -2\pi A\left(
t\right)  \right)  =-2\left(  L^{2}\left(  t\right)  -4\pi A\left(  t\right)
\right)  \leq0.
\]
Hence
\begin{equation}
0\leq L^{2}\left(  t\right)  -4\pi A\left(  t\right)  \leq\left(  L^{2}\left(
0\right)  -4\pi A\left(  0\right)  \right)  e^{-2t}\label{ipd}%
\end{equation}
as long as the solution exists.\ 

On the other hand if we compute the time derivative of the isoperimetric
ratio, we get
\begin{align*}
& \frac{d}{dt}\frac{L^{2}\left(  t\right)  }{4\pi A\left(  t\right)  }\\
& =\frac{4\pi A\left(  t\right)  \left(  2L^{2}\left(  t\right)  -4\pi
L\left(  t\right)  H\left(  L\left(  t\right)  ,A\left(  t\right)  \right)
\right)  -L^{2}\left(  t\right)  4\pi\left(  \int_{\gamma_{t}}\frac{1}%
{k}ds-H\left(  L\left(  t\right)  ,A\left(  t\right)  \right)  L\left(
t\right)  \right)  }{\left(  4\pi A\left(  t\right)  \right)  ^{2}}%
\end{align*}
and by the refined Green-Osher's inequality\ (see Lin-Tsai \cite{LT1})%
\begin{equation}
\int_{\gamma_{t}}\frac{1}{k}ds\geq\frac{2}{\pi}\left(  L^{2}\left(  t\right)
-4\pi A\left(  t\right)  \right)  +2A\left(  t\right)  \label{GO2}%
\end{equation}
we obtain%
\begin{equation}
\frac{d}{dt}\frac{L^{2}\left(  t\right)  }{4\pi A\left(  t\right)  }\leq
\frac{L\left(  t\right)  }{4\pi\left(  A\left(  t\right)  \right)  ^{2}%
}\left[  L^{2}\left(  t\right)  -4\pi A\left(  t\right)  \right]  \left[
H\left(  L\left(  t\right)  ,A\left(  t\right)  \right)  -\frac{2}{\pi
}L\left(  t\right)  \right]  .
\end{equation}
One can see that the three flows in (\ref{F4}),\ (\ref{F3}) and (\ref{F8}) are
all isoperimetric ratio decreasing.\ The isoperimetric ratio is also
decreasing as long as $H\ $is a negative function.\ See Remark \ref{rmk1}\ also.\ 

\begin{remark}
According to a private communication with S.-L.\ Pan,\ he also obtained the
same refined Green-Osher's inequality\ (\ref{GO2}) recently.\ Note that
(\ref{GO2})\ is an improvement of (\ref{GO1})\ by$\ \left(  L^{2}-4\pi
A\right)  /\pi\geq0$.\ Unlike the situation in most of the isoperimetric
inequalities where the equality cases occur only at circles, the equality case
of (\ref{GO2}) occurs if and only if the convex closed curve has support
function $u\left(  \theta\right)  $ given by%
\begin{equation}
u\left(  \theta\right)  =a_{0}+a_{1}\cos\theta+b_{1}\sin\theta+a_{2}%
\cos2\theta+b_{2}\sin2\theta,\ \ \ \theta\in\left[  0,2\pi\right]
\end{equation}
for some constants $a_{0},\ a_{1},\ b_{1},\ a_{2},\ b_{2}\ $satisfying%
\begin{equation}
u_{\theta\theta}\left(  \theta\right)  +u\left(  \theta\right)  =a_{0}%
-3a_{2}\cos2\theta-3b_{2}\sin2\theta>0\ \ \ \text{for all\ \ \ }\theta
\in\left[  0,2\pi\right]  .
\end{equation}
See Lin-Tsai \cite{LT1} for details.\ 
\end{remark}

\section{Convergence of the flow\ (\ref{main-1}).}

Assume that the solution $L\left(  t\right)  $ of the ODE (\ref{L-ODE}) is
defined on time interval$\ [0,\infty).\ $Then the function
\[
u\left(  \theta,t\right)  :=B\left(  \theta,t\right)  +\frac{L\left(
t\right)  }{2\pi}%
\]
is defined on $S^{1}\times\lbrack0,\infty)\ $and satisfies equation
(\ref{dudt}) everywhere. This$\ u\left(  \theta,t\right)  $ clearly satisfies
\begin{equation}
u_{\theta\theta}\left(  \theta,t\right)  +u\left(  \theta,t\right)
<\infty\ \ \ \text{for\ all\ \ \ }\left(  \theta,t\right)  \in S^{1}%
\times\lbrack0,\infty), \label{111}%
\end{equation}
which implies that
\begin{equation}
k\left(  \theta,t\right)  =\frac{1}{u_{\theta\theta}\left(  \theta,t\right)
+u\left(  \theta,t\right)  }>0\ \ \ \text{for\ all\ \ \ }\left(
\theta,t\right)  \in S^{1}\times\lbrack0,\infty). \label{333}%
\end{equation}
We also know that%
\begin{equation}
u_{\theta\theta}\left(  \theta,t\right)  +u\left(  \theta,t\right)
>0\ \ \ \text{for all\ \ \ }\left(  \theta,t\right)  \in S^{1}\times
\lbrack0,T_{\ast}) \label{222}%
\end{equation}
for some short time $T_{\ast}.\ $However, we can not exclude the possibility
that$\ u_{\theta\theta}\left(  \theta_{0},t_{0}\right)  +u\left(  \theta
_{0},t_{0}\right)  =0\ $at some$\ \left(  \theta_{0},t_{0}\right)  \in
S^{1}\times\lbrack0,\infty).\ $If this happens, then at time $t_{0}\ $we can
not use $u\left(  \theta,t_{0}\right)  $ to construct a smooth\ convex closed
curve $\gamma_{t_{0}}.$ In particular, this means that the flow (\ref{main-1})
may develop a singularity at time $t_{0}\ $with $k\left(  \theta_{0}%
,t_{0}\right)  =\infty\ $for some $\theta_{0}.\ $Note that
\begin{align}
&  u_{\theta\theta}\left(  \theta,t\right)  +u\left(  \theta,t\right)
\nonumber\\
&  =B_{\theta\theta}\left(  \theta,t\right)  +B\left(  \theta,t\right)
+\frac{L\left(  t\right)  }{2\pi}\nonumber\\
&  =\frac{e^{t}}{2\sqrt{\pi t}}\int_{-\infty}^{\infty}e^{-\frac{\xi^{2}}{4t}%
}\left(  u_{\theta\theta}\left(  \theta-\xi,0\right)  +u\left(  \theta
-\xi,0\right)  -\frac{1}{2\pi}\int_{0}^{2\pi}u\left(  \sigma,0\right)
d\sigma\right)  d\xi+\frac{L\left(  t\right)  }{2\pi} \label{GOOD}%
\end{align}
and,\ unless $\gamma_{0}$ is a circle, the function $u_{\theta\theta}\left(
\theta-\xi,0\right)  +u\left(  \theta-\xi,0\right)  -\frac{1}{2\pi}\int
_{0}^{2\pi}u\left(  \sigma,0\right)  d\sigma$ is somewhere positive and
somewhere negative,\ making it difficult to exclude the possibility.

We conclude the following:

\begin{theorem}
Let $H\left(  p,q\right)  :\left(  0,\infty\right)  \times\left(
0,\infty\right)  \rightarrow\mathbb{R}$ be a smooth function and let
$\gamma_{0}$ be a smooth convex closed curve.\ Assume that the ODE
(\ref{L-ODE})\ for $L\left(  t\right)  \ $is defined on $[0,\infty).$ Then the
nonlocal flow (\ref{main-1}) either develops a singularity (with $k=\infty
\ $somewhere) in finite time or the flow is defined on $S^{1}\times
\lbrack0,\infty),$ with each $\gamma_{t}\ $remaining smooth and convex,\ and
its support function\ $u\left(  \theta,t\right)  $ satisfies the following
$C^{\infty}\ $convergence on $S^{1}$:
\begin{align}
&  \lim_{t\rightarrow\infty}\left(  u\left(  \theta,t\right)  -\frac{L\left(
t\right)  }{2\pi}\right) \nonumber\\
&  =\left(  \frac{1}{\pi}\int_{0}^{2\pi}u\left(  \theta,0\right)  \cos\theta
d\theta\right)  \cos\theta+\left(  \frac{1}{\pi}\int_{0}^{2\pi}u\left(
\theta,0\right)  \sin\theta d\theta\right)  \sin\theta,\ \ \ \forall
\ \theta\in S^{1}. \label{converge}%
\end{align}

\end{theorem}

\begin{remark}
In the above theorem, the length $L\left(  t\right)  $ may go to infinity or
approach a positive\ constant or tend to zero as $t\rightarrow\infty
.\ $The$\ $same for the area $A\left(  t\right)  \ $due to (\ref{ipd}).$\ $The
geometric meaning of the above theorem is that $\gamma_{t}$ converges to a
circle $C_{t}\ $with length $L\left(  t\right)  \ $centered at the point
\begin{equation}
P=\left(  \frac{1}{\pi}\int_{0}^{2\pi}u\left(  \theta,0\right)  \cos\theta
d\theta,\ \frac{1}{\pi}\int_{0}^{2\pi}u\left(  \theta,0\right)  \sin\theta
d\theta\right)  \in\mathbb{R}^{2}. \label{P}%
\end{equation}

\end{remark}

\begin{remark}
\label{rmk1}If $H\left(  p,q\right)  :\left(  0,\infty\right)  \times\left(
0,\infty\right)  \rightarrow\left(  -\infty,0\right)  \ $is a negative
function, then the ODE (\ref{L-ODE}) implies that$\ L\left(  t\right)  \geq$
$L\left(  0\right)  e^{t}\ $for all $t\in\lbrack0,\infty).\ $Now
\begin{align*}
&  u_{\theta\theta}\left(  \theta,t\right)  +u\left(  \theta,t\right) \\
&  =\frac{e^{t}}{2\sqrt{\pi t}}\int_{-\infty}^{\infty}e^{-\frac{\xi^{2}}{4t}%
}\left(  u_{\theta\theta}\left(  \theta-\xi,0\right)  +u\left(  \theta
-\xi,0\right)  \right)  d\xi+\left(  \frac{L\left(  t\right)  }{2\pi}%
-\frac{L\left(  0\right)  }{2\pi}e^{t}\right)  >0
\end{align*}
due to $u_{\theta\theta}\left(  \theta,0\right)  +u\left(  \theta,0\right)
>0\ $everywhere.\ Hence it is impossible for the flow to develop a singularity
in finite time.\ 
\end{remark}

%

\proof
The proof is now straightforward.\ The convexity of $\gamma_{t}$ is due to
(\ref{333}).\ For the convergence, the equation for $u\left(  \theta,t\right)
-L\left(  t\right)  /2\pi\ $is linear with
\[
\int_{0}^{2\pi}\left(  u\left(  \theta,t\right)  -\frac{L\left(  t\right)
}{2\pi}\right)  d\theta=0\ \ \ \text{for all\ \ \ }t\in\lbrack0,\infty).
\]
Using Fourier series expansion for$\ u\left(  \theta,t\right)  -L\left(
t\right)  /2\pi,\ $the convergence result follows from standard theory. $%
\hfill
\square$

\ \ \ \ 

For the ODE (\ref{L-ODE}), it is also possible that $L\left(  t\right)  $ is
defined only on a finite time interval $[0,T_{\max})$ with either
$\lim_{t\rightarrow T_{\max}}L\left(  t\right)  =\infty\ $or $\lim
_{t\rightarrow T_{\max}}L\left(  t\right)  =0.$ In the first case, we note
that%
\begin{equation}
\left\vert u\left(  \theta,t\right)  -\frac{L\left(  t\right)  }{2\pi
}\right\vert \leq Ce^{t}\ \ \ \text{for\ all\ \ \ }t\in\lbrack0,T_{\max}),
\end{equation}
where $C>0$ is a constant depending only on the initial
curve.\ Hence$\ u\left(  \theta,t\right)  \rightarrow\infty$ uniformly on
$S^{1}\ $as$\ t\rightarrow T_{\max}.\ $By%
\begin{equation}
\frac{\partial}{\partial t}\left(  \frac{\partial^{m}u}{\partial\theta^{m}%
}\right)  \left(  \theta,t\right)  =\left(  \frac{\partial^{m}u}%
{\partial\theta^{m}}\right)  _{\theta\theta}\left(  \theta,t\right)  +\left(
\frac{\partial^{m}u}{\partial\theta^{m}}\right)  \left(  \theta,t\right)
,\ \ \ \int_{0}^{2\pi}\frac{\partial^{m}u}{\partial\theta^{m}}\left(
\theta,t\right)  d\theta=0
\end{equation}
and using Fourier series expansion (or by the result in
Chow-Gulliver\ \cite{CG})\ we know that for each$\ m\in\mathbb{N},\ \left\vert
\left(  \partial^{m}u/\partial\theta^{m}\right)  \left(  \theta,t\right)
\right\vert $ is uniformly bounded on $S^{1}\times\lbrack0,T_{\max}).\ $Hence
if we rescale the curve $\hat{\gamma}_{t}$ by considering$\ \hat{\gamma}%
_{t}=2\pi\gamma_{t}/L\left(  t\right)  ,\ $its support function $\hat
{u}\left(  \theta,t\right)  $ will satisfy$\ \hat{u}\left(  \theta,t\right)
\rightarrow1\ $in$\ C^{\infty}(S^{1})\ $as$\ t\rightarrow T_{\max}.\ $

We conclude the following:

\begin{theorem}
If the ODE (\ref{L-ODE})\ for $L\left(  t\right)  \ $is defined on
$[0,T_{\max})\ $with$\ T_{\max}<\infty\ $and$\ \lim_{t\rightarrow T_{\max}%
}L\left(  t\right)  =\infty.$ Then the nonlocal flow (\ref{main-1}) either
develops a singularity\ before $T_{\max}\ $or the flow is defined on
$S^{1}\times\lbrack0,T_{\max}),$ with each $\gamma_{t}\ $remaining smooth and
convex,\ and the support function\ $\hat{u}\left(  \theta,t\right)  $ of the
rescaled curve $\hat{\gamma}_{t}=2\pi\gamma_{t}/L\left(  t\right)
\ $satisfies$\ \hat{u}\left(  \theta,t\right)  \rightarrow1\ $in$\ C^{\infty
}(S^{1})\ $as$\ t\rightarrow T_{\max}.\ $
\end{theorem}

For the case $\lim_{t\rightarrow T_{\max}}L\left(  t\right)  =0,$ by
(\ref{int-eq}) we have%
\[
\frac{1}{k\left(  \theta,t\right)  }=Z\left(  \theta,t\right)  +\frac{L\left(
t\right)  }{2\pi},\ \ \ Z\left(  \theta,t\right)  =e^{t}\int_{-\infty}%
^{\infty}\frac{1}{2\sqrt{\pi t}}e^{-\frac{\left(  \theta-\xi\right)  ^{2}}%
{4t}}\left(  \frac{1}{k\left(  \xi,0\right)  }-\frac{1}{2\pi}\int_{0}^{2\pi
}\frac{1}{k\left(  \sigma,0\right)  }d\sigma\right)  d\xi
\]
where the function\ $Z\left(  \theta,t\right)  $ is somewhere$\ $positive and
somewhere negative (unless $\gamma_{0}$ is a circle)$\ $for each $t\in
\lbrack0,\infty)\ $due to $\int_{0}^{2\pi}Z\left(  \theta,t\right)
d\theta=0.\ $Therefore before time $T_{\max}$ we must have $k=\infty
\ $somewhere.\ Hence we conclude:

\begin{theorem}
If the ODE (\ref{L-ODE})\ for $L\left(  t\right)  \ $is defined on
$[0,T_{\max})\ $with$\ T_{\max}<\infty\ $and$\ \lim_{t\rightarrow T_{\max}%
}L\left(  t\right)  =0.$ Then if $\gamma_{0}$ is not a circle,\ the nonlocal
flow (\ref{main-1}) must develop a singularity\ before time $T_{\max}\ $with
$k=\infty\ $somewhere.
\end{theorem}

If the ODE (\ref{L-ODE})\ for $L\left(  t\right)  \ $is defined only on a
maximal domain$\ [0,T_{\max})\ $with$\ T_{\max}<\infty,\ $then $\left(
L\left(  t\right)  ,A\left(  t\right)  \right)  $ must approach the boundary
of the domain $\left(  0,\infty\right)  \times\left(  0,\infty\right)  \ $of
$H\left(  p,q\right)  .\ $If$\ \lim_{t\rightarrow T_{\max}}L\left(  t\right)
=\infty\ $(or $0$)$,\ $then the same for$\ \lim_{t\rightarrow T_{\max}%
}A\left(  t\right)  $ due to (\ref{ipd})\ and the isoperimetric inequality.
There remains the last case that$\ \lim_{t\rightarrow T_{\max}}L\left(
t\right)  =\ell\in\left(  0,\infty\right)  \ $but$\ \lim_{t\rightarrow
T_{\max}}A\left(  t\right)  =0.\ $In such a case, by Gage's inequality for
convex closed curves%
\[
\int_{\gamma_{t}}k^{2}ds\geq\frac{\pi L\left(  t\right)  }{A\left(  t\right)
}\rightarrow\infty\ \ \ \text{as\ \ \ }t\rightarrow T_{\max}%
\]
we must have$\ \lim_{t\rightarrow T_{\max}}k_{\max}\left(  t\right)
=\infty.\ $We conclude:

\begin{theorem}
Assume the ODE (\ref{L-ODE})\ for $L\left(  t\right)  \ $is
defined\ on\ $[0,T_{\max})\ $with$\ T_{\max}<\infty\ $and$\ \lim_{t\rightarrow
T_{\max}}L\left(  t\right)  =\ell\in\left(  0,\infty\right)  .\ $%
If$\ \lim_{t\rightarrow T_{\max}}A\left(  t\right)  =0,\ $then we must
have$\ \lim_{t\rightarrow T_{\max}}k_{\max}\left(  t\right)  =\infty.\ $
\end{theorem}

\section{Conclusion.}

We can allow the nonlocal term $H\left(  L\left(  t\right)  ,A\left(
t\right)  \right)  $ to be more general. For example in Ma-Cheng\ \cite{MC}%
,\ they considered an interesting nonlocal flow of the form%
\begin{equation}
\dfrac{\partial X}{\partial t}\left(  \varphi,t\right)  =\left(  \frac
{1}{L\left(  t\right)  }\int_{0}^{L\left(  t\right)  }\frac{1}{k\left(
\varphi,t\right)  }ds-\dfrac{1}{k\left(  \varphi,t\right)  }\right)
\mathbf{N}_{in}\left(  \varphi,t\right)  , \label{MC-flow}%
\end{equation}
which is an \emph{area-preserving} flow.\ For this flow the support function
$u\left(  \theta,t\right)  \ $and the length $L\left(  t\right)  $ still
satisfies the same\ equation (\ref{u-L}), and so$\ u\left(  \theta,t\right)
=L\left(  t\right)  /2\pi+B\left(  \theta,t\right)  ,\ $where$\ B\left(
\theta,t\right)  $ is given by (\ref{uL}). The ODE for $L\left(  t\right)  $
in this case is
\begin{align}
\frac{dL}{dt}\left(  t\right)   &  =L\left(  t\right)  -\frac{2\pi}{L\left(
t\right)  }\int_{0}^{L\left(  t\right)  }\frac{1}{k\left(  \varphi,t\right)
}ds=L\left(  t\right)  -\frac{2\pi}{L\left(  t\right)  }\int_{0}^{2\theta
}\left[  u_{\theta\theta}\left(  \theta,t\right)  +u\left(  \theta,t\right)
\right]  ^{2}d\theta\nonumber\\
&  =L\left(  t\right)  -\frac{2\pi}{L\left(  t\right)  }\int_{0}^{2\theta
}\left(  B_{\theta\theta}\left(  \theta,t\right)  +B\left(  \theta,t\right)
+\frac{L\left(  t\right)  }{2\pi}\right)  ^{2}d\theta. \label{MC}%
\end{align}
Again (\ref{MC}) is a self-contained ODE and has short time existence.\ Thus
Ma-Cheng's flow has solution for short time.\ In view of this, we see that the
nonlocal flow (\ref{main-1}) has short time existence as long as the
quantities in the nonlocal term $H\ $can be expressed in terms of the support
function $u\left(  \theta,t\right)  .\ $Since the position vector of a convex
closed curve $\gamma_{t}$ is uniquely determined by its support function via
the formula (\ref{position}), we can conclude that essentially \emph{any}
$1/k$-type nonlocal\ flow has short time existence and can be dealt with by
the above ODE method.\ 

\bigskip

\ \ \ \ \ \ 

\textbf{Acknowledgments.\ \ }\ \ Both authors are grateful for the support of
NSC (National Science Council)\ of Taiwan, under grant numbers
99-2115-M-007-001\ and 99-2115-M-007-MY3. We also appreciate the support of
NCTS\ (National Center of Theoretical Sciences)\ of Taiwan during the past
several years.\ The writing of this paper is motivated by the two papers by
Pan-Yang\ \cite{PY} and Ma-Cheng\ \cite{MC}, we thank these authors.\ 

\ \ \ 

\ \ \ \ \ \

\ \ \ \ 

\ \ \ \ \ \ 

Yu-Chu Lin

Department of Mathematics

National Tsing Hua University

Hsinchu 30013,\ TAIWAN

E-mail:\ \textit{yclin@math.nthu.edu.tw}

\ 

\ \ \ 

\ \ \ \ \ 

Dong-Ho Tsai\ 

Department of Mathematics

National Tsing Hua University

Hsinchu 30013,\ TAIWAN

E-mail:\ \textit{dhtsai@math.nthu.edu.tw}

\bigskip\

\end{document}